\documentclass[12pt]{article}
\usepackage{amssymb}

\def\header{A.Miller\hfill A Dedekind Finite Borel Set \hfill}
 \def\header{\today} 

\def\Si{\Sigma}

\def\al{\alpha}
\def\bb{{\mathcal B}}
\def\be{\beta}
\def\cl{{\rm cl}}
\def\concat#1{{\hat{\phantom{a}}\la #1\ra}}
\def\concatx{{\hat{\phantom{a}}}}
\def\de{\delta}
\def\emp{\emptyset}
\def\ff{{\mathcal F}}
\def\fix{{\rm fix}}
\def\forces{{\;\Vdash}}

\def\gg{{\mathcal G}}
\def\hh{{\mathcal H}}
\def\ka{\kappa}
\def\lam{\lambda}
\def\la{\langle}

\def\name#1{\stackrel{\circ}{#1}}

\def\nn{{\mathcal N}}
\def\om{\omega}
\def\one{{\bf 1}}
\def\poset{{\mathbb P}}
\def\proof{\par\noindent Proof\par\noindent}
\def\pr{\prime}
\def\qed{\par\noindent QED\par\bigskip}

\def\ra{\rangle}
\def\ra{\rangle}
\def\res{\upharpoonright}
\def\rmand{\mbox{ and }}
\def\rmiff{{\mbox{ iff }}}
\def\rmor{\mbox{ or }}
\def\rr{{\mathcal R}}
\def\si{\sigma}
\def\sm{{\setminus}}
\def\st{\;:\;} 
\def\sumposet{\Si_{n<\om}\poset(D+E_n)}
\def\su{\subseteq}

\newtheorem{theorem}{Theorem}[section]

\newtheorem{prop}[theorem]{Proposition}
\newtheorem{cor}[theorem]{Corollary}
\newtheorem{define}[theorem]{Definition}
\newtheorem{lemma}[theorem]{Lemma}
\newtheorem{remark}[theorem]{Remark}

\begin{document}

\begin{center}
{\large A Dedekind Finite Borel Set}
\end{center}

\begin{flushright}
Arnold W. Miller
\footnote{
\par We would like to thank Mark Fuller and Daniel Kane 
for bringing this problem to
our attention and especially Mark for his enthusiastic interest.
}
\end{flushright}

\def\address{\begin{flushleft}
Arnold W. Miller \\
miller@math.wisc.edu \\
http://www.math.wisc.edu/$\sim$miller\\
University of Wisconsin-Madison \\
Department of Mathematics, Van Vleck Hall \\
480 Lincoln Drive \\
Madison, Wisconsin 53706-1388 \\
\end{flushleft}}

\begin{center}  Abstract  \end{center}
\begin{quote}
In this paper we prove three theorems about the theory of
Borel sets in models of ZF without any form of the axiom of choice.
We prove that if $B\su 2^\om$ is a $G_{\de\si}$-set then either
$B$ is countable or $B$ contains a perfect subset.  Second, we prove
that if $2^\om$ is the countable union of countable sets, then
there exists an $F_{\si\de}$ set $C\su 2^\om$ such that $C$ is
uncountable but contains no perfect subset.  Finally, we construct
a model of ZF in which we have an
infinite Dedekind finite $D\su 2^\om$ 
which is $F_{\si\de}$.
\end{quote}

{\small \tableofcontents}

\section{ Introduction }

In this paper we assume the theory ZF but we do not assume any form of
the axiom of choice, in particular, we do not assume the countable
axiom of choice (which says that choice functions exist for countable
families of nonempty sets). For example, we do not assume that the
countable union of countable sets is countable.

It is well-known that assuming the countable axiom of choice that
every uncountable Borel set contains a perfect set.  In fact, it is
not hard to see, that assuming the countable axiom of choice that
every Borel subset of $2^\om$ is the projection of a closed subset of
$2^\om\times \om^\om$, i.e., an analytic set, and that every
uncountable analytic set contains a perfect set.

\begin{define}
\begin{enumerate}
\item For $s\in 2^{<\om}$ define the basic clopen set:
$$[s]=\{x\in 2^\om\st s\su x\}.$$
\item A set $U\su 2^\om$ is open iff it is the union of
basic clopen sets.
\item A set $A\su 2^\om$ is $G_{\de}$ iff it is the intersection of
a countable family of open sets.
\item A set $B\su 2^\om$ is $G_{\de\si}$ iff it is the union of
a countable family of $G_\de$-sets.
\item Similarly define $F$ be the closed sets, i.e., complements of
open sets, $F_\si$ the countable
unions of closed sets, and $F_{\si\de}$ the 
countable intersections of $F_\si$'s.
\item A subset $P\su 2^\om$ is perfect iff it is homeomorphic to $2^\om$.
\end{enumerate}
\end{define}

\begin{theorem}\label{perfthm}
If $A\su 2^\om$ is a $G_{\de\si}$ set, then $A$ is
countable or contains a perfect set.
\end{theorem}

\begin{theorem} \label{thmctble}
Suppose that $2^\om$ is the countable union of countable sets.
Then there exists an $F_{\si\de}$ set $B\su 2^\om$ which is uncountable
but contains no perfect subset.
\end{theorem}

In the Feferman-Levy model the $2^\om$ is the countable union
of countable sets (see Cohen \cite{cohen} p.143, Jech
\cite{jech} p.142).
Note that this implies that every set
$B\su 2^\om$ is the countable union of countable sets.  Since
a countable subset of $2^\om$ is an $F_\si$, it follows that
every subset of $2^\om$ is $F_{\si\si}$, i.e., a countable
union of countable
unions of closed sets.  By taking complements every subset of $2^\om$
is $G_{\de\de}$.  So the set $B$ in Theorem \ref{thmctble} is
$F_{\si\de}$, $F_{\si\si}$, and $G_{\de\de}$.

In ZF without using any choice at all there
exists a $G_{\de\si}$-set which is not $F_{\si\de}$, see
Theorem 2.1 of Miller \cite{longbor}.

A set $D$ is Dedekind finite iff every one-one map of $D$ into itself
is onto.  Equivalently, there is no one-one map of $\om$ into $D$.
Assuming the axiom of choice every Dedekind finite set is finite.
The book Herrlich \cite{herrlich} pp.43-50 summarizes many 
of the basic results about Dedekind finite sets.

By  infinite set we simply mean that the set is
not finite, i.e.,
cannot be put into one-to-one correspondence with some
finite ordinal $n\in\om$.

\begin{theorem}\label{thmded}
Suppose that $M$ is a countable transitive model of ZF and
$$M\models D\su 2^\om \mbox{ is an infinite Dedekind finite set. }$$
Then there exists a symmetric submodel $\nn$ of a generic extension
of $M$ such that
$$\nn\models D \mbox{ is a Dedekind finite $F_{\si\de}$-set.}$$
\end{theorem}

For our forcing terminology over models of ZF see Miller
\cite{longbor} section 3.

\begin{remark}
If $2^\om$ is the countable union of countable sets, then
there are no infinite Dedekind finite $D\su 2^\om$.
This is because the countable union of finite subsets of
a linearly orderable set is countable.
\end{remark}

Besides the notion of Dedekind finite there are many other
``definitions of finiteness'', i.e.,
properties which are equivalent to finite assuming the
axiom of choice
(see Truss \cite{trussded}, L\'{e}vy \cite{levy},
Howard and Yorke \cite{howard}, De la Cruz \cite{cruz}).
Most of  them are
inconsistent with being an infinite subset of $2^\om$.  One exception
is is $\Delta_5$  (see Truss
\cite{trussded}):

 A set $D$ is $\Delta_5$ iff there does not exist an
onto map $f:D\to D\cup\{*\}$ where $*$ is not an element of $D$.  

It
is possible to have an infinite $\Delta_5$ subset of $2^\om$.
Let us say $D\su 2^\om$ has the density-Dedekind property iff it is a
dense subset of $2^\om$ and for any $E\su D$ there
exists an open set $U\su 2^\om$ such that $d\in E$ iff $d\in U\cap D$
for all but finitely many $d\in D$.  Density-Dedekind
implies $\Delta_5$.  In the basic Cohen model of ZF in which 
choice fails
(see Jech \cite{jech} p.66-68) there is a generic
Dedekind finite set $A\su 2^\om$.
It is not hard to show that in fact $A$ 
has the density-Dedekind property and hence
is $\Delta_5$.  The notion of density-Dedekind seems to us to
be analogous to that of Luzin set in set theory with choice.

We don't know if it is possible to have an infinite Borel
$\Delta_5$-set.  Almost-disjoint sets forcing destroys the
density-Dedekind property.

A set is amorphous iff every subset of it is finite or cofinite.
This is analogous in model theory with the Baldwin and Lachlan
notion of strongly minimal set (see Truss \cite{trussamor}, Creed,
Truss \cite{creedomin}, Mendick, Truss \cite{mend}, and
Walczak-Typke \cite{wal}).
An infinite $D\su 2^\om$ cannot be amorphous.  We don't know if there
could be an uncountable Borel set $D\su 2^\om$ such that every
subset is countable or co-countable (i.e., quasi-amorphous, see
Creed, Truss \cite{creedquasi}).  

Monro \cite{monro} constructed Dedekind finite sets which are large
in the sense that they can be mapped onto a cardinal $\ka$.
The ones he constructed were subsets of $2^\ka$.
It is possible to have a Dedekind finite Borel set which maps
onto $\om_1$ (or any other larger $\om_\al$ if desired). 
By Theorem \ref{thmded} it is enough to find a
Dedekind finite set $D\su 2^\om$ which maps onto $\om_1$.
Such a $D$ can be constructed by using
a slight variant of the second Cohen model, see Jech
\cite{jech} pp. 68-71.
 
In computability theory, the notion of Dedekind finite
is analogous
to that of Dekker's notion of an isol.  There are over 180 of papers
on the theory of isols, although currently the subject seems to have
fallen out of fashion.   Two which connect the theory of isols and
Dedekind finite cardinals are Ellentuck \cite{ellen} and McCarty
\cite{mccarty}.  Perhaps there are analogies between Borel
Dedekind finite sets and co-simple isols,
i.e., complements of simple sets. See for example, Downey and
Slaman \cite{downey} which contains work on co-simple isols.

\section{Proof of Theorem \ref{perfthm}}

\begin{define} Recall the following: \label{treedef}
\begin{enumerate}
\item A nonempty $T\su 2^{<\om}$ is a tree iff 
$\forall s,t\in 2^{<\om}$ if $s\su t\in T$, then $s\in T$.
\item For $T$ a tree
$$[T]=\{x\in 2^\om\st \forall n<\om\;\; x\res n\in T\}$$
\item For $T$ a tree and $s\in T$
$$T(s)=\{t\in T\st t\su s\rmor s\su t\}.$$
\item $T$ is  perfect iff 
$\forall s\in T\;\;\exists t\in T \;\;s\su t$
and both $t\concat 0\in T$ and 
$t\concat 1\in T$.
\end{enumerate}
\end{define}

\noindent The proof of the following proposition
is left to the reader.

\begin{prop}
A set $C\su 2^\om$ is closed iff there exists a tree
$T\su 2^{<\om}$ such that $C=[T]$.
A set $P\su 2^\om$ is perfect iff there is a perfect tree
$T\su 2^{<\om}$ such that $P=[T]$.
In both cases we may demand that the tree $T$ have no terminal
nodes, i.e., for any $s\in T$ either $s\concat 0\in T$ or
$s\concat 1\in T$.
\end{prop}

\begin{lemma}\label{perflem1}
Let $\bb$ be the family
of nonempty countable closed subsets of $2^\om$.
Then there is a function $\ff:\bb\to (2^\om)^\om$ such that
if $\ff(C)=f$, then $f:\om\to C$ is an onto map.
\end{lemma}
\proof

This argument is ancient set theory, the
Cantor-Bendixson derivative.
(Recall we must not use of the axiom of choice.)

Let $C$ be a nonempty countable closed set. Define
$$T=\{s\in 2^{<\om}\st [s]\cap C\neq\emp\}.$$
Hence $[T]=C$.

Inductively define a sequence of trees $T_\al\su 2^{<\om}$ 
for $\al$ an ordinal as follows:
\begin{enumerate}
\item $T_0=T$
\item $T_\lam=\bigcap_{\al<\lam} T_\al$ is $\lam$ is a limit ordinal
\item $T_{\al+1}=T_\al\sm \{t\in T_\al\st |[T_\al(t)]| \leq 1\}.$
\end{enumerate}
Note that $\al\leq \be$ implies $T_\be\su T_\al$.

If $T_{\al+1}=T_\al$, then $T_\al=T_\be$ for all
$\be>\al$.  By the replacement axiom there must be an ordinal
$\al$ such that $T_{\al+1}=T_\al$.   Since $T_\al\su T$ we
have that $[T_\al]\su C$ and since $C$ is countable, it must
be that $T_\al$ is empty, since otherwise it is easy
to check that it is a perfect tree.

For each $x\in C$ there exists a unique ordinal $\al_x<\al$ such
that 
$$x\in [T_{\al_x}]\sm [T_{\al_x+1}].$$
Let $n$ be the least such that $x\res n \notin T_{\al_x+1}$ and
put $s_x=x\res n$.  We claim that the map $q:C\to 2^{<\om}$ defined
by $q(x)=s_x$ is one-to-one.  To see this suppose that
$s_x=s_y$.  If $\al_x<\al_y$, then we get a contradiction since
$s_x\notin T_{\al_x+1}$ and $T_{\al_y}\su T_{\al_x+1}$.  So
$\al_x=\al_y$ and from the definition of $T_{\al_x+1}$ we see
that $x=y$.

To get our onto map $f:\om\to C$, let $x_0$ be the lexicographically
least element of $C$ and let $\{t_n\st n<\om\}$ be a fixed enumeration of
$2^{<\om}$.  Given any $n$ if $t_n=s_x$ for some $x\in C$ let
$f(n)=x$ and otherwise let $f(n)=x_0$.

No choice is being used in our definition of $f$, so we may
define $\ff(C)=f$.
\qed

\begin{cor}
The countable union of closed
subsets of $2^\om$ each of which is countable is countable.
\end{cor}

\begin{lemma}\label{perflem2}
Let $\hh$ be the family of nonempty countable
$G_\de$ subsets of $2^\om$.  Then there is a
function $\gg:\hh\to (2^\om)^\om$ such that
if
$\gg(H)=g$, then $g:\om\to H$ is an onto map.
\end{lemma}
\proof
This argument is also ancient set theory (although perhaps not
as well known), the Hausdorff difference hierarchy.  Hausdorff
proved that disjoint $G_\de$ sets can be separated by a set
which is in the difference hierarchy of closed sets
(see Kechris \cite{kechris} p.176).

Let $H,K\su 2^\om$ be disjoint $G_\de$-sets.
Define closed sets $C_\al\su 2^\om$ for $\al$ an ordinal as follows:

 $C_0=\cl(H)$  (we use $\cl(X)$ to denote the closure of $X$)

 $C_1=\cl(C_0\cap K)$
 
 $C_2=\cl(C_1\cap H)$

 $\vdots$
 
 $C_\om=\bigcap_{n<\om}C_n$,
 
 $C_{\om+1}=cl(C_\om\cap K)$

 $C_{\om+2}=cl(C_\om\cap H)$

and so forth, in general, for $\lambda$ a limit ordinal and $n<\om$:

$$C_\lam=\bigcap_{\al<\lam}C_\al$$
$$C_{\lam+n+1}=\left\{
\begin{array}{ll}
\cl(C_{\lam+n}\cap H) & \mbox{if $n$ is odd}\\
\cl(C_{\lam+n}\cap K) & \mbox{if $n$ is even}\\
\end{array}\right.$$

It is clear that if $\al\leq \be$ then $C_\be\su C_\al$.   Also
if\footnote{With
a little more work it is enough that $C_\al=C_{\al+1}$.}
 $C_\al=C_{\al+2}$ then for all $\be>\al$, $C_\be=C_\al$.
Hence there must be an ordinal $\al_0$ such that
$C_{\al_0}=C_\be$ for all $\be>\al_0$.

We claim that $C_{\al_0}$ is empty,  otherwise,
$H$ and $K$ are both dense in it.  Hence it would follow 
that $G\cap H\neq \emp$.
To see this let $$T=\{s\in 2^{<\om}\st [s]\cap C_{\al_0}\neq\emp\}.$$
Write $H=\bigcap_{n<\om} U_n$ and $K=\bigcap_{n<\om} V_n$ where
$U_n$ and $V_n$ are open sets.  Since $H$ and $K$ are dense in
$C_\al$, it must be that for every $s\in T$ and $n<\om$,
there exists $t\in T$
with $s\su t$ and $[t]\su U_n\cap V_n$.   But now it is easy to
construct $x\in [T]\cap H\cap K$.

Since $C_{\al_0}$ is empty we have that the difference sets:
$$D=\bigcup\{(C_\al\sm C_{\al+1}) \st \al \mbox{ is even }\}$$
and
$$E=(2^\om\sm C_0)\cup\bigcup \{(C_\al\sm C_{\al+1}) 
\st \al \mbox{ is odd }\}.$$
are complementary\footnote{The ordinal $\al_0$ must 
be countable and the unions could
be taken over $\al\leq \al_0+2$ but we don't need this for
our proof.}.
We claim that $H\su D$ and $K\su E$.  To see why, suppose that
$x\in H$.  Since $C_0=\cl(H)$ it must be that there is some
ordinal $\al$ such that $x\in C_\al\sm C_{\al+1}$.  This $\al$ cannot
be odd, since $C_{\al+1}=\cl(C_\al\cap H)$.

Now suppose that $H$ is a countable $G_\de$-set.  Then
 $K=2^\om\sm H$ is also a $G_\de$-set.  From which it
follows that

$$H\;\;= \;\;\bigcup \{(C_\al\sm C_{\al+1}) \st \al
\mbox{ is even }\}.$$

Define  
$$T_\al=\{s\in 2^{<\om}\st [s]\cap C_\al\neq\emp\}.$$
So for $\al<\al_0$ we have that each $T_\al$ is a nonempty tree
without terminal nodes such that $C_\al=[T_\al]$.
For $s\in 2^{<\om}$ with length greater than $0$,
let $s^*\su s$ with $|s^*|=|s|-1$. Let
$$Q_\al=\{s\in T_\al\sm T_{\al+1}\st s^*\in T_{\al+1}\}$$
i.e, the minimal nodes of $T_\al\sm T_{\al+1}$.

For each even $\al$ since $C_\al\sm C_{\al+1}\su H$ and $H$ is
countable we have that $[T_\al(s)]$ is a countable set
for each $s\in Q_\al$. Note that the $Q_\al$ are pairwise disjoint.  Let
$Q\su 2^{<\om}$ be the set of all $s$ such that $s\in Q_\al$ and
$\al$ is even.  For each $s\in Q$ define $f_s:\om\to 2^{<\om}$ by
$\ff([T_\al(s)])=f_s$ where $s\in Q_\al$.  It follows that
the map $h:Q\times\om\to H$ defined by $h(s,n)=f_s(n)$ is onto
$H$ and may easily be readjusted to an onto map $g:\om\to H$.
Put $\gg(H)=g$.
\qed

\begin{cor}\label{perfcor}
The countable union of countable $G_\de$ subsets of $2^\om$ is countable.
\end{cor}
\proof
Suppose that $\bigcup_{n<\om}H_n$ is given
where each $H_n$ is a countable $G_\de$-set.  Let
$\gg(H_n)=g_n$.  Then define an onto map 
$$g:\om\times\om\to \bigcup_{n<\om}H_n
\;\;\;\;\mbox{ by }\;\;\;\; g(n,m)=g_n(m).$$
\qed

\bigskip

Proof of Theorem \ref{perfthm}.

\bigskip

It follows immediately from Corollary \ref{perfcor} that
we need only show that an uncountable $G_\de$-set $H\su 2^\om$
must contain a perfect set.

Define
$$H^\pr=\{x\in H\st \forall n<\om\;\;\; ([x\res n]\cap H) \mbox { is uncountable }
\}.$$
Note that $H^\pr$ is nonempty, since otherwise
$$H = \bigcup\{[s]\cap H\st [s]\cap H \mbox{ is countable},\;\; s\in 2^{<\om}\}$$
and since any set of the form $[x\res n]\cap H$ is $G_\de$ and
the countable union of $G_\de$-sets is countable, we would get a contradiction.

Define
$$T=\{s\in 2^{<\om}\st [s]\cap H^\pr\neq \emp\}.$$
We claim that $T$ is a perfect tree.  To see this suppose
that $s\in T$. Then $s$ will have incompatible extensions in $T$ unless
$H^\pr\cap [s]=\{x\}$.
This would mean that for every extension $t$ of $s$ which is incomparable
to $x$ that $H\cap [t]$ is countable.
But since $s\su x$ we know that
$H\cap [s]$ is uncountable.  But this
contradicts the fact that the countable union of countable
$G_\de$-sets is countable.

Now suppose $H=\bigcap_{n<\om} U_n$ where each $U_n$ is open.  We
construct $$(s_\si\in T:\si\in 2^{<\om})$$ by induction on the
length of $\si$.  Given $s_\si$ with $|\si|=n$ let $t\in T(s_\si)$ be
the first in some fixed ordering of $2^{<\om}$ with $[t]\su U_n$.  
Then using that $T$ is perfect similarly find
$s_{\si\concat i}$ for $i=0,1$ incomparable extensions of $t$.
Then 
$$T^\pr=\{t\st \exists\si\;\; t\su s_\si\}$$
is a perfect subtree of $T$ such that $[T^\pr]\su H$.
\qed

\section{Proof of Theorem \ref{thmctble}}

\begin{define}
Let $\la,\ra:\om\times\om\to\om$ be a fixed bijection,
i.e., a pairing function.  For each $n\in\om$
define the map $\pi_n:2^\om\to 2^\om$ by:
$$\pi_n(x)=y \mbox { iff } \forall m\in\om\;\;y(m)=x(\la n,m\ra).$$
\end{define}

\begin{lemma}\label{lemnice}
Suppose that $2^\om$ is the countable union of countable sets. 
Then there exists $(F_n\st n\in\om)$ such that
\begin{enumerate}
\item $2^\om=\bigcup_{n<\om}F_n$ and each $F_n$ is countable,
\item $F_n$ is a proper subset of $F_{n+1}$ for each $n$, and
\item $F_n$ is closed under $\pi_m$ for all $n,m<\om$.
\end{enumerate}
\end{lemma}
\proof
Define a map $H:\om^{<\om}\times 2^{\om}\to 2^{\om}$ inductively
by $$H(\la \ra,x)=x \rmand
 H(s\concat m,x)=\pi_m(H(s,x)).$$  Given that $2^\om=\bigcup_{n<\om}L_n$
where each $L_n$ is countable, let
$$F_n=H(\om^{<\om}\times \bigcup_{m\leq n}L_m).$$
Then the $F_n$ are countable, increasing, cover $2^\om$, and
closed under the projection maps $\pi_m$.  To get them to be
properly increasing just pass to a subsequence.
\qed

Define
$$B_n=\{x\in 2^\om:\pi_n(x)\in F_{n+1}\sm F_n
\rmor [\pi_n(x)\in F_n \rmand \pi_n(x)=\pi_{n+1}(x)]
\}.$$
Note that each $B_n$ is an $F_\si$-set.
Let $B=\bigcap_{n<\om} B_n$.

The set $B$ is uncountable because
there is a map $h$ from $B$ onto $2^\om$.  Define
$h$ by $h(x)=\pi_n(x)$ iff $\pi_n(x)=\pi_m(x)$ for all $m>n$.  Such an $n$ must
exists because for any $x$ there exists $n$ such that
$x\in F_n$ and hence $\pi_m(x)\in F_n$ for all $m$.  It is easy to check
that $h$ maps $B$ onto $2^\om$.

But $B$ cannot contain a perfect set.  Suppose for
contradiction that $T\su 2^{<\om}$
is a perfect tree and $[T]\su B$.
For each $x\in [T]$ define $h(x)=n$ to be the least $n$ so that 
$\pi_n(x)=\pi_m(x)$ for all $m>n$.  For any $n$ the set
of all $x\in [T]$ with $h(x)\leq n$ is closed.  By Corollary \ref{perfcor}
it must be that for some $n$ that there exists a perfect subtree
$T^\pr\su T$ such that $x\in [T^\pr]$ implies $h(x)< n$.  But
the map
$$k:[T^\pr]\to \prod_{m<n}F_m \mbox{ defined by } k(x) = (\pi_m(x):m<n)$$
would map a perfect set one-one into a countable set.
\qed

\begin{remark}
We don't really need Corollary \ref{perfcor} in the above
proof, since it is easy to show that a perfect set cannot be the
countable union of countable closed sets.  For example, each would
have to be nowhere dense.
\end{remark}

\begin{remark}
In the Feferman-Levy model the set $B$ has the stronger
property that there is no one-one map (continuous or not)
taking $2^\om$ into $B$.  Also Lemma \ref{lemnice} is trivially
true in that model since we take $F_n=M[G_n]\cap 2^\om$.
\end{remark}

\section{Proof of Theorem \ref{thmded}}

\begin{define}
A poset $\poset$ is $\si$-centered iff there exists
$(\Si_n:n<\om)$ such that $\poset=\bigcup_{n<\om}\Si_n$
and each $\Si_n$ is centered, i.e., for any finite $F\su\Si_n$ there
exists $p\in\poset$ such that $p\leq q$ for every $q\in F$.
\end{define}

We begin with a preservation lemma:

\begin{lemma}\label{sigmacent}
Suppose that $M$ is a countable transitive model of ZF and
$$M\models \mbox{$\poset$ is $\si$-centered and $D$ 
is Dedekind finite. }$$
Then for any $G$ $\poset$-generic over $M$
$$M[G]\models \mbox{$D$ 
is Dedekind finite. }$$
\end{lemma}
\proof
Working in $M$ let $(\Si_n:n<\om)$ witness the $\si$-centeredness
of $\poset$.
Suppose for contradiction that
$$p\forces \name{f}:\check{\om}\to \check{D}\mbox{ is one-one.}$$

Define
$$D_{n,m}=\{x\in D\st 
\exists q\in \Si_n\;\; q\leq p \rmand q\forces \name{f}(m)=\check{x}
\}.$$
Since $\Si_n$ is centered, $|D_{n,m}|\leq 1$.  Since
$D$ is Dedekind finite, the set $$E=\bigcup_{n,m<\om}D_{n,m}$$ is
finite.  But
$$p\forces \mbox{ the range of $\name{f}$ is a subset of $\check{E}$}$$
which is a contradiction.
\qed

Remark.  To preserve the Dedekind finiteness of $D\su 2^\om$ 
it would be enough to assume that $\poset=\bigcup_{n<\om}\Si_n$ where
each $\Si_n$ had the n-c.c., i.e., no antichain of size greater than
$n$.

Next we give a description of the well-known almost-disjoint sets 
forcing of Solovay.

\begin{define}\label{defposet}
For $A\su 2^\om$ define
$$\poset(A)=\{\la Q,F\ra\st Q\su 2^{<\om},\; F\su A, \rmand
\mbox{ both $Q$ and $F$ are finite }\}.$$
For $p,q\in \poset(A)$ define
$p\leq q$ iff $Q_q\su Q_p$, $F_q\su F_p$, and
$s\not\su x$ for all $s\in Q_p\sm Q_q$ and $x\in F_q$.
\end{define}

We use $\one=(\emp,\emp)$ to denote the trivial element of
$\poset(A)$.

Note that $\la Q,F_1\ra\leq \la Q,F_2\ra$ whenever $F_2\su F_1$.
Hence given $Q\su 2^{<\om}$ finite, if we define 
$$\Si_Q=\{p\in\poset(A)\st Q_p=Q\}$$
then $\Si_Q$ is centered and
$$\poset(A)=\bigcup\{\Si_Q\st Q\su 2^{<\om}\mbox{ is finite }\}$$ 
shows that $\poset(A)$ is $\si$-centered.

If $G$ is $\poset(A)$-generic over $M$, then we can define
$$R=R^G=\bigcup\{Q_p\st p\in G\}.$$
Easy density arguments show that for every $x\in 2^\om\cap M$
\begin{itemize}
\item if $x\in A$, then $\{n\st x\res n\in R\}$ is finite, and
\item if $x\notin A$, then $\{n\st x\res n\in R\}$ is infinite.
\end{itemize}

Next we consider automorphisms of the poset $\poset(A)$.

\begin{define}
A map $\hat{\pi}:2^{<\om}\to 2^{<\om}$ is a tree automorphism iff
$\hat{\pi}$ is a bijection such that for all $s,t\in 2^{<\om}$
$$s\su t \rmiff \hat{\pi}(s)\su \hat{\pi}(t).$$ 
\end{define}

A tree automorphism $\hat{\pi}$ induces a map from $2^\om$ to itself
by letting $\hat{\pi}(x)=y$ where $y$ is determined by 
$\hat{\pi}(x\res n)=y\res n$ for every $n<\om$. 

\begin{lemma}
Suppose $\hat{\pi}$ is a tree automorphism such that 
$\hat{\pi}(x)\in A$
for every $x\in A$.  Then $\pi:\poset(A)\to \poset(A)$ defined
by 
$$\pi(Q,F)=(\{\hat{\pi}(s):s\in Q\},\{\hat{\pi}(x):x\in F\})$$
is an automorphism of $\poset(A)$.
\end{lemma}
\proof
We need to show that
$$p\leq q \rmiff \pi(p)\leq\pi(q).$$
It is easy to check that
$$Q_q\su Q_p \rmiff \hat{\pi}(Q_q)\su \hat{\pi}(Q_p)$$
and
$$F_q\su F_p \rmiff \hat{\pi}(F_q)\su \hat{\pi}(F_p).$$
For the third clause in definition \ref{defposet} 
note that for $s\in 2^{<\om}$ and
$x\in 2^\om$ that
$$s\su x  \rmiff \hat{\pi}(s)\su \hat{\pi}(x).$$

\qed

\begin{define}
For any $n<\om$ define 
$$E_n=\{x\in 2^\om\st \forall k <n \;\; x(k)=0
\rmand \exists l>n\; \forall k>l \;\; x(k)=0\}$$
As usual for $x,y\in 2^\om$ define $x+y$ to be their pointwise sum mod 2, i.e.,
$$ \forall n\;\;\; (x+y)(n)\equiv x(n)+y(n) \mbox{ mod }2$$
and for $A,B\su 2^\om$ define
$$A+B=\{x+y\st x\in A \rmand y\in B\}.$$
\end{define}

\begin{lemma}
For $D\su 2^\om$ Dedekind finite
$$D=\bigcap_{n<\om}(D+E_n).$$
\end{lemma}
\proof
Since the constant zero function is in every $E_n$ it is clear
that 
$$D\su\bigcap_{n<\om}(D+E_n).$$

Now suppose for contradiction that $x\in \bigcap_{n<\om}(D+E_n)$ but
$x\notin D$.  
Consider the equivalence class of $x$ under ``equal mod finite'': 
$x+E_0$.
Since this class can be well-ordered in type $\om$ we know
that the set: $$F=D\cap (x+E_0)$$ is finite.  Take $n<\om$ large enough so
that for all $u,v\in F\cup\{x\}$ if $u\res n= v\res n$ then $u=v$.  
But $x\in D+E_n$ which means that there exists $d\in D$ with
$d\res n= x\res n$.   But $d\in F$ which is a contradiction.

\qed

\begin{define}
We define the poset $\poset$ to be the direct sum of the posets:
$\poset(D+E_n)$, i.e.,
$$\poset=\sumposet.$$
This means $p\in \poset$ iff $p=(p_n:n<\om)$ where each $p_n\in \poset(D+E_n)$
and $p_n=\one$ for all but finitely many $n$.  It is ordered
coordinatewise: 
$$p\leq q \rmiff p_n\leq q_n \mbox{ for all }n.$$
As before, given any $G$ a $\sumposet$-filter and $n<\om$ we define
$$R_n=R_n^G=\{s\in 2^{<\om}\st \exists p\in G \mbox{ with } s\in Q_{p_n}\}.$$
\end{define}

It is clear that for $G$ a $\poset$-generic filter over $M$ that
for every $n$ and $x\in D+E_n$ there are at most finitely many $k<\om$
with $x\res k\in R_n$.

\begin{lemma}
The poset $\poset=\sumposet$ is $\si$-centered.
\end{lemma}
\proof
For any finite sequence $\vec{W}=(Q_i:i<n)$ of finite subsets of $2^{<\om}$
define
$$\Si_{\vec{Q}}=\{p\in\poset\st \forall i<n\;\; Q_{p_i}=Q_i \;\rmand\;
\forall i\geq n\;\;
p_i=\one\}.$$
Then each $\Si_{\vec{Q}}$ is centered and $\poset$ is the countable
union of them.
\qed

\begin{define}
For $R\su 2^{<\om}$ define
$$H(R)=\{x\in 2^\om\st \exists^\infty k\;\; x\res k\in R\}$$
Here $\exists^\infty k$ stands for ``there exists infinitely many $k$''.
\end{define}

\begin{lemma}\label{lemhr}
For $R\su 2^{<\om}$ the set $H(R)$ is a $G_\de$-set.
Suppose $\rr$ is a countable family of subsets of $2^{<\om}$, then
$\bigcap\{H(R)\st R\in\rr\}$
is a $G_\de$-set.
\end{lemma}
\proof
$$H(R)=\bigcap_{n<\om}\bigcup\{[s] \st s\in R\rmand |s|>n\}.$$
Letting $\rr=\{R_n\st n<\om\}$ we have that
$$\bigcap\{H(R)\st R\in\rr\}=\bigcap_{n,m<\om} 
\bigcup\{[s] \st s\in R_n\rmand |s|>m\}.$$
\qed

Note that $H(R_n)$ is a $G_\de$-set disjoint from $D+E_n$.  Our goal
is to make the complement of $D$ to be a countable union of $G_\de$ sets
in a symmetric submodel of $M[G]$.

We describe the automorphisms of $\poset$ which we will use.
\begin{define}
\begin{enumerate}
\item For $s\in 2^{<\om}$ define $\hat{\pi}_s:2^{<\om}\to 2^{<\om}$ to
be the tree automorphism which swaps $s\concat 0$ and
$s\concat 1$, i.e., 
$$\hat{\pi}(r)=\left\{
\begin{array}{ll}
s \concat {1-i} \concatx t & \mbox{ if } r=s \concat i \concatx t\\
r & \mbox{ if $r$ does not extend $s$. }\\
\end{array}
\right.
$$
\item For each $n$ we let $\gg_n$ be the group of automorphisms of
$\poset(D+E_n)$ which are generated by
$\{\pi_s\st s\in 2^{<\om}\rmand |s|>n\}$.
\item We take $\gg$ to be the direct sum of the $\gg_n$, i.e.,
$\pi\in\gg$ iff $\pi=(\pi_n:n<\om)$ where each $\pi_n\in\gg_n$ and
$\pi_n$ is the identity except for finitely many $n$.
\item We take $\ff$ to be the filter of subgroups of $\gg$ which is generated
by $\{H_n:n<\om\}$ where 
$$H_n=\{\pi\in\gg\st \forall m<n \;\;\pi_m \mbox{ is the identity }\}.$$
\end{enumerate}
\end{define}

It is easy to check that $\ff$ is a normal filter.

We use the terminology $\hat{\pi}$ (a hatted $\pi$) 
to denote tree automorphisms and unhatted $\pi$'s to
denote the corresponding automorphism of $\poset$ and the action on the
$\poset$-names. 
We use $\nn$ to denote the symmetric model
$M\su \nn\su M[G]$.
We use the terminology $\fix(\tau)$ to denote
the subgroup of $\gg$ which fixes the $\poset$-name $\tau$.

Let 
$$\name{R}_n=\{(p,\check{s})\st p\in\poset\rmand s\in Q_{p_n}\}$$
then $H_{n+1}\su \fix(\name{R}_n)$ and so $R_n\in \nn$.
The following lemma is key:

\begin{lemma}\label{keylem}
Given $p\in\poset$, $\name{x}$, and $n_0<\om$ such that
$H_{n_0}\su \fix(\name{x})$ and
$$p\forces \name{x}\in 2^\om\sm (D+E_{n_0})$$
then
$$p\forces \exists^\infty k\;\; \name{x}\res k\in\name{R}_{n_0}.$$
\end{lemma}
\proof
If not there exists $q\leq p$ and $N>n_0$ such that
$$q\forces \forall n>\check{N}\;\; \name{x}\res n\notin \name{R}_{n_0}.$$

\bigskip

\noindent{\bf Claim}.  There exists $r\leq q$ and $s,t_0,t_1\in 2^{<\om}$ with
\begin{enumerate}
\item $|s|>N>n_0$
\item $\{t\in 2^{<\om}\st s\su t\}\cap Q_{q_{n_0}}=\emp$
\item $[s]\cap F_{q_{n_0}}=\emp$.
\item $s\su t_0$, $s\su t_1$,  $|t_0|= |t_1|$
\item $t_0\in Q_{r_{n_0}}$
\item $r\forces \check{t}_1\su\name{x}$
\end{enumerate}

\bigskip
Since $p$ is forcing that $x$ is not in $D+E_{n_0}$ it easy to find
$r_1\leq q$ and $s$ such that
$$r_1\forces \check{s}\su\name{x}$$
and $s$ satisfies 1,2, and 3.
Next choose any $t_0$ with $s\su t_0$ and
$t_0\not\su y$ for all $y\in F_{r_{1,n_0}}$ and put
$r_2=r_1$ except 
$$Q_{r_{2,n_0}}=Q_{r_{1,n_0}}\cup\{t_0\}.$$
Finally
find  $r\leq r_2$ and $t_1$ with $|t_0|=t_1$ and
$r\forces \check{t}_1\su\name{x}$. This proves the Claim.

\bigskip

  Now find $\hat{\pi}$ a tree automorphism
in $\hat{\gg}_{n_0}$ such 
$\hat{\pi}(t_0)=t_1$ and 
fixes all $t$ except for possibly those extending
$s$.  A precise description would be to let:
$$\{n\st t_0(n)\neq t_1(n)\}=\{n_1<n_2<\ldots<n_k\}$$
then 
$$\hat{\pi}=\hat{\pi}_{s_k}\circ\hat{\pi}_{s_{k-1}}\circ\cdots\circ
\hat{\pi}_{s_1}$$ 
where $s_i=t_1\res n_i$.  Note that $\pi\in\gg_{n_0}$ because
$|t_0|=|t_1|\geq |s|>N>n_0$ so necessarily $n_1>n_0$.  
Let $\pi\in\gg$ also name 
the automorphism of $\poset$ which is $\pi$ on the $n_0^{th}$
coordinate and the identity on all other coordinates.  Then
$\pi\in H_{n_0}$ and hence $\pi(\name{x})=\name{x}$ and
so by (6) of the Claim
$$\pi(r)\forces \check{t}_1\su \name{x}.$$
Note that by (2) and (3) of the Claim, we have $\pi(q)=q$ and so $\pi(r)\leq q$
and thus:
$$\pi(r)\forces\forall n>N\;\; \name{x}\res n\notin \name{R}_{n_0}.$$
By (5) of the Claim and the definition of $\pi$ we
have that $t_1\in Q_{\pi(r)_{n_0}}$ so we have:
$$\pi(r)\forces \check{t}_1\in  \name{R}_{n_0}.$$
But $|t_1|>N$ gives us a contradiction.
\qed

Let
$$\rr_n=\{\hat{\pi}(R_n)\st \hat{\pi}\in\hat{\gg}_{n}\}.$$
That is we take the set of all images of $R_n$ under the tree 
automorphisms which determine
$\gg_n$.   Since each $R_n$ is in $\nn$ and 
$\hat{\gg}_{n}$ is in the ground model, it is clear that
each $\rr_n$ is in $\nn$.

\begin{lemma}\label{lemravel}
For each $\pi\in \gg$ and $n<\om$:
$$\pi(\name{R}_n)^{G}=\hat{\pi}_n^{-1}(R_n).$$
\end{lemma}
\proof
This amounts to unraveling the definitions.  The following are
equivalent:
\begin{itemize}
\item $s\in \pi(\name{R}_n)^{G}$
\item $\exists p\in G$ such that $(p,s)\in \pi(\name{R}_n)$ and
 $(p,s)=(\pi(q),s)$ where $s\in Q_{q_n}$
\item $\exists p\in G$ such that $s\in Q_{\pi_{n}^{-1}(p)}$
(equivalently $\pi_{n}(s)\in Q_{p_n}$)
\item $\pi_n(s)\in R_n$
\item $s\in \hat{\pi}_n^{-1}(R_n)$.
\end{itemize}
\qed

\begin{lemma}
The sequence $(\rr_n:n\in\om)$ is in $\nn$.
\end{lemma}
\proof
Letting
$$\name{\rr}_n=\{\pi(\name{R}_n)\st \pi\in \gg\}$$
we see that $\fix(\name{\rr}_n)=\gg$ for every $n$,
hence the $\om$-sequence has a name fixed by every $\pi$ in $\gg$.
\qed

Next we show that in the hypothesis of the key lemma (Lemma \ref{keylem})
we may assume that the
trivial condition $\one$ is doing the forcing.

\begin{lemma} \label{lemone}
Fix $G$ a $\poset$-filter generic over $M$.
Suppose $x\in (2^\om\sm D)\cap\nn$.  Then $x$ has a
hereditarily symmetric
name $\name{x}$ for which there is an $n_0$ such that
$H_{n_0}\su\fix(\name{x})$ and
$$\one\forces \name{x}\in 2^\om\sm (D+E_{n_0}).$$
\end{lemma}
\proof
Let $\tau$ be any hereditarily symmetric name for $x$, i.e.,
$\tau^G=x$.  Let $p\in G$ and $n_0$ be such that $H_{n_0}\su\fix(\tau)$
and
$$p\forces \tau\in 2^\om\sm (D+E_{n_0}).$$

Now work in the ground model $M$.
Fix $z\in M\cap (2^\om\sm (D+E_{n_0}))$.
In $M$ define $\name{x}$ to the set of all $(q,\la m,i\ra)$
such that either
$$q\forces \mbox{``} \tau\in 2^\om\sm (D+E_{n_0}) \wedge \tau(m)=i
\mbox{''}$$
or
$$z(m)=i \rmand q\forces \neg (\tau\in 2^\om\sm (D+E_{n_0})).$$

For our particular $G$, since $p\in G$ the second clause is never
invoked when evaluating $\name{x}^G$, hence
$x=\tau^G=\name{x}^G$.  Clearly, $H_{n_0}\su\fix(\tau)\su\fix(\name{x})$.
Finally $\one$ forces what it should because for any generic
filter $G^\pr$ either
the first clause is invoked and $\tau^{G^\pr}=\name{x}^{G^\pr}$
or the second clause is invoked and
$\name{x}^{G^\pr}=z$ where $z$ was chosen to be in $2^\om\sm (D+E_{n_0})$.
\qed

\begin{lemma}
$\nn\models (2^\om\sm D)=\bigcup_{n<\om}\bigcap_{R\in \rr_n}H(R)$.
\end{lemma}
\proof
Recall that 
$H(R_n)$ is a $G_\de$-set which is disjoint 
from $D+E_n$ and
hence from $D$.
For any $R\in\rr_n$ we have that $R=\hat{\pi}(R_n)$ for
some $\pi\in\gg_n$.  But by Lemma \ref{lemravel}
$$R=\hat{\pi}(R_n)=\pi(\name{R}_n)^{G}=\name{R_n}^{\pi^{-1}(G)}$$
and so
$H(R)$ is disjoint from $D$.

Conversely suppose in $\nn$ that $x\in (2^\om\sm D)$.  Then
by Lemma \ref{lemone}
$x$ has a name $\name{x}$ for which  there exists $n_0$ such that
$H_{n_0}\su\fix(\name{x})$ and
$$\one\forces \name{x}\in 2^\om\sm (D+E_{n_0}).$$
By the key Lemma \ref{keylem}
$$\one\forces \exists^\infty k\;\; \name{x}\res k\in\name{R}_{n_0}$$
i.e.,
$$\one\forces \name{x}\in H(\name{R}_{n_0}).$$
Since $H_{n_0}\su\fix(\name{x})$ we have that
$$\one\forces \name{x}\in H(\pi(\name{R}_{n_0}))$$
for all $\pi\in H_{n_0}$ and so it follows from Lemma \ref{lemravel} 
that $$x\in \bigcap_{R\in \rr_n}H(R).$$
\qed

It follows from this Lemma that in $\nn$ the complement of $D$ is
a $G_{\de\si}$ set and hence $D$ is an $F_{\si\de}$-set.
Since  $\poset$ is $\si$-centered
we have that 
$$M[G]\models D\mbox{ is Dedekind finite}$$
and since $M\su\nn\su M[G]$
$$\nn\models D\mbox{ is a Dedekind finite $F_{\si\de}$-set. }$$

This concludes the proof of Theorem \ref{thmded}.

\address


\addcontentsline{toc}{section}{References}

\begin{thebibliography}{99}

\bibitem{cohen} Cohen, Paul J.; {\bf Set theory and the
continuum hypothesis.} W. A. Benjamin, Inc., 
New York-Amsterdam 1966 vi+154 pp.

\bibitem{creedomin} Creed, P.; Truss, J. K.; On o-amorphous sets. Ann. Pure
Appl. Logic 101 (2000), no. 2-3, 185--226.

\bibitem{creedquasi} Creed, P.; Truss, J. K.; On quasi-amorphous sets. Arch.
Math. Logic 40 (2001), no. 8, 581--596.

\bibitem{cruz} De la Cruz, Omar; Finiteness and choice. Fund. Math. 173 (2002),
no. 1, 57--76. 

\bibitem{downey} Downey, Rod; Slaman, Theodore A.;
On co-simple isols and their intersection types.
Ann. Pure Appl. Logic 56 (1992), no. 1-3, 221--237.

\bibitem{ellen} Ellentuck, Erik; The universal properties of Dedekind finite
cardinals. Ann. of Math. (2) 82 1965 225--248.

\bibitem{herrlich} Herrlich, Horst; {\bf Axiom of choice}. Lecture Notes in
Mathematics, 1876. Springer-Verlag, Berlin, 2006. xiv+194 pp. ISBN:
978-3-540-30989-5

\bibitem{howard} Howard, Paul E.; Yorke, Mary F.; Definitions of finite. Fund.
Math. 133 (1989), no. 3, 169--177.

\bibitem{jech} Jech, Thomas J.; {\bf The axiom of choice.} Studies in Logic and
the Foundations of Mathematics, Vol. 75. North-Holland Publishing Co.,
Amsterdam-London; Amercan Elsevier Publishing Co., Inc., New York, 1973. xi+202
pp.

\bibitem{kechris} Kechris, Alexander S.; {\bf Classical descriptive set
theory.} Graduate Texts in Mathematics, 156. Springer-Verlag,
New York, 1995. xviii+402 pp.

\bibitem{levy} L\'{e}vy, A.; The independence of various definitions of
finiteness. Fund. Math. 46 1958 1--13.

\bibitem{mend} Mendick, G. S.; Truss, J. K.; A notion of rank in set theory
without choice. Arch. Math. Logic 42 (2003), no. 2, 165--178.

\bibitem{longbor} Miller, Arnold W.;  Long Borel hierarchies, Math Logic
Quarterly, 54(2008), 301-316.

\bibitem{mccarty} McCarty, Charles;
Markov's principle, isols and Dedekind finite sets.
J. Symbolic Logic 53 (1988), no. 4, 1042--1069.

\bibitem{monro} Monro, G. P.; Independence results concerning Dedekind-finite
sets. J. Austral. Math. Soc. 19 (1975), 35--46.

\bibitem{trussded} Truss, J. K.; Classes of Dedekind finite cardinals. Fund. Math.
84 (1974), no. 3, 187--208.

\bibitem{trussamor} Truss, J. K.; The structure of amorphous sets.  Ann. Pure
Appl. Logic  73  (1995),  no. 2, 191--233.

\bibitem{wal} Walczak-Typke, A. C.; The first-order structure of weakly
Dedekind-finite sets.  J. Symbolic Logic  70  (2005),  no. 4, 1161--1170.

\end{thebibliography}
\end{document}